\documentclass[twoside,11pt]{amsart}
\usepackage{amsmath,latexsym,amssymb}
\usepackage{emlines}

\newcommand{\rar}{\rightarrow}
\newcommand{\lar}{\longrightarrow}

\newtheorem{Theorem}{Theorem}[section]

\newtheorem{Corollary}[Theorem]{Corollary}
\newtheorem{Example}[Theorem]{Example}
\newtheorem{Proposition}[Theorem]{Proposition}

\newtheorem{Question}[Theorem]{Question}
\newtheorem{Remark}[Theorem]{Remark}

\newtheorem{Definition}[Theorem]{Definition}

\def\demo{\noindent{\bf Proof. }}
\def\QED{\hfill$\Box$}

\def\red{\mbox{\rm r}}

\def\depth{\mbox{\rm depth }}
\def\height{\mbox{\rm height }}

\newcommand{\Rees}{\mbox{${\mathcal R}$}}
\def\gr{\mbox{\rm gr}}

\newcommand{\xvec}[1]{\ensuremath{x_{1},\ldots,x_{#1}}}
\def\m{{\mathfrak m}}
\def\b{{\mathfrak b}}
\def\a{{\mathfrak a}}
\def\bar#1{{\overline{#1}}}
 
\begin{document}

\title[Indices of Normalization of Ideals]
{\Large\bf Indices of normalization of ideals}

\author[C. Polini, B. Ulrich,  W.V. Vasconcelos and R. Villarreal]
{C. Polini \and B. Ulrich \and W. V. Vasconcelos
\and R. Villarreal}

\thanks{AMS 2000 {\em Mathematics Subject Classification}.
Primary 13A30; Secondary 13B22, 13H10, 13H15.}

\thanks{The first author 
 gratefully acknowledge partial
support from the NSA, the second and third authors from the NSF,
 and the last author from SNI, Mexico. }

\address{Department of Mathematics, University of Notre Dame,
Notre Dame, Indiana 46556} \email{cpolini@nd.edu}
\urladdr{www.nd.edu/{\textasciitilde}cpolini}

\address{Department of Mathematics, Purdue University,
West Lafayette, Indiana 47907} \email{ulrich@math.purdue.edu}
\urladdr{www.math.purdue.edu/{\textasciitilde}ulrich}

\address{Department of Mathematics, Rutgers University,
Piscataway, New Jersey 08854} \email{vasconce@math.rutgers.edu}
\urladdr{www.math.rutgers.edu/{\textasciitilde}vasconce}

\address{Departamento de
Matem\'aticas,
 Centro de Investigaci\'on y de Estudios
Avanzados del
IPN,
Apartado Postal
14--740,
 07000 Mexico City, D.F.}
\email{vila@math.cinvestav.mx}\urladdr{www.math.cinvestav.mx/{\textasciitilde}vila}

\vspace{0.1in}

\begin{abstract} We derive numerical estimates controlling
the  intertwined properties of the normalization of an ideal and of
the computational complexity of general processes for its
construction. In \cite{ni1}, this goal was carried out for
equimultiple ideals via the examination of Hilbert functions.
Here we add to this picture, in an important case, 
  how certain Hilbert functions provide a description of
the locations of the generators of the normalization of ideals of
dimension zero. 
 We also present a rare
instance of normalization of a class of homogeneous ideals by a
single colon operation.
\end{abstract}

\maketitle

\smallskip

\section{Introduction}

Let $R$ be a Noetherian integral domain and let $I$ be an ideal. The
{\em normalization} of $I$ is the integral closure in $R[t]$,
$\overline{A}$, 
of the Rees algebra
$A=R[It]$ of $I$. In case $R$ is normal the nuance disappears.
The properties of $\overline{A}$ add significantly
to an understanding of $I$ and of the constructions it supports. The
{\em index} terminology refers to the integers related to the
description and construction of
\[ \bar{A}= \sum_{n\geq 0} \overline{I^n}t^n =
{R}[\overline{I}t, \ldots,
\overline{I^{s_0}}t^{s_0}]
.\]

In addition to the overall task of describing the generators and
relations of $\bar{A}$, it includes the understanding of the
following quantities:
\begin{itemize}
\item[{\rm (i)}] 
Numerical indices for equalities of the type: 
find $s$ such that $$(\bar{A})_{n+s}= (A)_n\cdot (\bar{A})_s, \quad n\geq
0.$$
\item[{\rm (ii)}]
Estimation of the number of  steps that effective processes must
traverse   between $A$ and
$\bar{A}$, 
\[A= A_0 \subset A_1 \subset \cdots \subset A_{r-1}\subset
A_r =
\bar{A}.\]
\item[{\rm (iii)}]
Express $r$, $s$ and $s_0$ in
terms of invariants of $A$.
\item[{\rm (iv)}]
Generators of $\bar{A}$: number of generators  and
distribution of their degrees in cases of interest.
\end{itemize}

\medskip

These general questions acquire a high degree of specificity when
$A=R[It]$, and the goal becomes the estimation of these indices in
terms of invariants of $I$. A general treatment of item (ii) is given
in \cite{jdeg1} and \cite{jdeg2}. For monomial ideals of finite
colength a  discussion is carried out in \cite{hilbsam}.
This paper is a sequel to \cite{ni1}, where some of the notions
developed here originated. The focus in \cite{ni1} was on
deriving bounds on the coefficient $e_1(\overline{A})$ of the Hilbert
function associated to ideals of finite co-length in local rings,
and its utilization in the estimation of the length of general
normalization algorithms. Here we introduce complementary notions and
use them to address some of the same goals for more general ideals,
but also show how known initial knowledge about the normalization
allows  us to give fairly detailed description of $\overline{A}$, 
 particularly
those affecting the
distribution of  its generators.

\smallskip

We now describe  the organization of the paper. Section \ref{normid} gives the
precise definitions of the indices outlined above (with one exception
best left for Section \ref{sallynormid}), and describe some relationships amongst
them. These indices acquire a sharp relief when the normalization
$\sum_{n\geq 0} \overline{I^n}t^n $ is Cohen-Macaulay
(Theorem~\ref{redintclos}). This result, whose proof 
follows {\it ipsis literis} the characterization of
Cohen-Macaulayness in the Rees algebras of $I$-adic filtrations
(\cite{AHT}, \cite{JK95}, \cite{SUV2}), has
various consequences.      
It is partly used to motivate the treatment in Section
\ref{sallynormid} of the Sally module of
the normalization algebra as a vehicle to study the number of
generators and their degrees. In case the associated graded ring of the
integral closure filtration $\mathcal{F}$ of an
$\mathfrak{m}$-primary ideal $I$, $\gr_{\mathcal{F}}(R)$ has depth
at least $\dim R-1$, there are several positivity relations on
the Hilbert coefficients, leading to  descriptions of the
distribution of the new generators (usually fewer as the degrees go
up), and overall bounds for their numbers.

In  Section \ref{onestepnorm},  we present one of the rare instances where the
normalization of the blowup ring is computed
using an explicit expression as a colon ideal. Our formula applies
to homogeneous ideals that are generated by forms of the same degree
and satisfy some additional assumptions.

\section{Normalization of Ideals}\label{normid}
This section introduces auxiliary constructions and
devices to examine   the
integral closure of ideals, and  to study the
properties and applications to normal ideals.

\subsection*{Indices of normalization}

We begin by introducing some measures for the normalization of
ideals.
Suppose $R$ is a commutative ring and $J,I$ are ideals of $R$ with
$J\subset I$. $J$ is a {\em reduction} of $I$ if $I^{r+1}=JI^r$ for
some integer $r$; the least such integer is the {\em reduction
number} of $I$ relative to $J$. It is denoted
$\red_J(I)$. $I$ is {\em equimultiple} if there is a reduction $J$
generated by $\height I$ elements.

\medskip

\begin{Definition}{\rm Let $R$ be a locally analytically unramified normal domain and
let $I$ be an ideal.
\begin{itemize}
\item[{\rm (i)}]
 The {\em  normalization
index}\index{normalization of an ideal!normalization index} of $I$ is the
smallest integer $s=s(I)$
such that
\[\overline{I^{n+1}} = I\cdot \overline{I^n} \quad  n\geq s.\]
\item[{\rm (ii)}]
 The {\em generation
index}\index{normalization of an ideal!generation index} of $I$ is the
smallest integer $s_0=s_0(I)$ such that
\[ \sum_{n\geq 0}\overline{I^n}t^n=
R[\overline{I}t, \ldots,
\overline{I^{s_0}}t^{s_0}].\]
\item[{\rm (iii)}] The {\em normal relation
 type}\index{normalization of an ideal!normal relation type} of $I$ is
 the  maximum degree of a minimal generating set of  the presentation ideal
\[ 0 \rar M \lar R[T_1, \ldots, T_m] \lar
R[\overline{I}t, \ldots,
\overline{I^{s_0}}t^{s_0}] \rar 0.\]
\end{itemize}
}\end{Definition}

\medskip

For example, if $R = k[x_1, \ldots, x_d]$ is a polynomial ring over
a field and $I = (x_1^d, \ldots, x_d^d)$, then $I_1 = \overline{I} =
(x_1, \ldots, x_d)^d$. It follows that $s_0(I) = 1$, while $s(I) =
\red_I(I_1) = d-1$.

\medskip

If $(R, \m)$ is a local ring, these indices have an expression in
term of the special fiber ring $F$ of the normalization map $A
\rightarrow \overline{A}$.

\smallskip

\begin{Proposition} \label{indicesnorm} With the above assumptions let \[ F =
\overline{A}/(\mathfrak{m}, It)\overline{A}= \sum_{n\geq 0} F_n.\]
We have
\begin{eqnarray*}
s(I) & = & \sup\{ n \mid F_n \neq 0\}, \\
s_0(I) & = & \inf \{n \mid F = F_0[F_1, \ldots, F_n]\}.
\end{eqnarray*}
Furthermore, if the index of nilpotency of $F_i$ is $r_i$, then
\[ s(I) \leq \sum_{i=1}^{s_0(I)} (r_i-1).\]
\end{Proposition}

\medskip

Although these integers are well defined---since $\overline{A}$ is
finite over $A$---it is not clear, even in case $R$ is a regular
local ring, which invariants of $R$ and of $I$ have a bearing on the
determination of $s(I)$. An affirmative case is that of a monomial
ideal $I$ of a ring of polynomials in $d$ indeterminates over a
field---when $s\leq d-1$ (according to
Corollary~\ref{normalpowers}).

\subsection*{Equimultiple ideals}  For primary ideals and some other
equimultiple ideals there are relations between the two indices of
normalization.

\begin{Proposition} Let $(R, \mathfrak{m})$ be a
local analytically unramified normal Cohen--Macaulay ring such that
 $\mathfrak{m}$ is a normal ideal. Let $I$ be
$\mathfrak{m}$--primary ideal with multiplicity $e(I)$. Then
\[ s(I)\leq e(I) ((s_0(I)+1)^d-1)-s_0(I)(2d-1).\]
\end{Proposition}

\demo Without loss of generality, we may assume that the residue
field of $R$ is infinite. Following Proposition~\ref{indicesnorm},
we estimate $s(I)$ (the Castelnuovo--Mumford regularity of $F$) in
terms of the indices of nilpotency of the components $F_n$, for
$n\leq s_0(I)$.

Let $J= (z_1, \ldots, z_d)$ be a minimal reduction of $I$. For each
component $I_n= \overline{I^n}$ of $\overline{A}$, we collect the
following data:
\begin{eqnarray*}
J_n & = & (z_1^n, \ldots, z_d^n), \  \textrm{a minimal reduction
of $I_n$}\\
e(I_n) & = & e(I)n^d, \   \textrm{the multiplicity of $I_n$} \\
r_n=r_{J_n}(I_n) & \leq & \frac{e(I_n)}{n}d - 2d + 1, \  \textrm {a
bound on the reduction number of $I_n$}.
\end{eqnarray*}
The last assertion follows from \cite[Theorem~7.14]{bookthree},
 once it is
observed that $I_n\subset \overline{\mathfrak{m}^n}= \mathfrak{m}^n$,
by the normality of $\mathfrak{m}$.

We are now ready to estimate the index of nilpotency of the component
$F_n$.
With the notation above, we have ${I_n}^{r_n+1}= J_nI_n^{r_n}$. When
 this relation is read  in $F$, it means that
$r_n+1\geq \textrm{index of nilpotency of $F_n$}$.

Following Proposition~\ref{indicesnorm}, we have
\[ s(I) \leq \sum_{n=1}^{s_0(I)} r_n = \sum_{n=1}^{s_0(I)} e(I) d
n^{d-1} -s_0(I)(2d-1),
\]
which we approximate with an elementary integral to get the
assertion.
 \QED

\bigskip

We can do considerably better when $R$ is a ring of polynomials
 over a field of characteristic zero.

\begin{Theorem} Let $R= k[x_1, \ldots, x_d]$ be a polynomial ring over a field of
characteristic zero and let $I$ be a homogeneous ideal that is
$(x_1, \ldots, x_d)$--primary. One has
\[ s(I) \leq (e(I)-1)s_0(I).\]
\end{Theorem}

\demo
 We begin by localizing
$R$ at the maximal homogeneous ideal and picking a minimal reduction
$J$ of $I$. We denote the associated graded ring of the filtration of
integral closures $\{ I_n= \overline{I^n}\}$ by $G$,
\[ G = \sum_{n\geq 0} I_n/I_{n+1}.\]
In this affine ring we can take for a Noether normalization a ring
$A = k[z_1, \ldots, z_d]$, where the $z_i$'s are the images in $G_1$
of a minimal set of generators of $J$.

There are two basic algebraic facts about the algebra $G$. First,
its multiplicity as a graded $A$--module is the same as that of the
associated graded ring of $I$, that is, $e(I)$. Second, since the
Rees algebra of the integral closure filtration is a normal domain,
so is the extended Rees algebra
\[ C=\sum_{n\in \mathbb{Z}}I_{n}t^n,\]
where we set $I_n=R$ for $n \leq 0$. Consequently the algebra $G
=C/(t^{-1})$ will satisfy the condition $S_1$ of Serre. This means
that as a module over $A$, $C$ is torsionfree.

We now apply the theory of Cayley-Hamilton equations to the elements
of the components of $G$ (see \cite[Chapter 9]{compu}):  For
 $u\in G_n$, we have an equation of integrality over $A$
\[ u^r + a_1u^{r-1} + \cdots + a_r=0,\]
where $a_i$ are homogeneous forms of $A$, in particular $a_i\in
A_{ni}$,
 and $r\leq e(G)= e(I)$.
Since $k$ has characteristic zero, using the
argument of
 \cite[Proposition 9.3.5]{compu}, we obtain an equality
\[ G_n^r = A_nG_n^{r-1}.\]
At the level of the filtration, this equality means that
\[ I_n^r \subset J^nI_{n}^{r-1} + I_{nr+1},\]
 which we weaken by
\[ I_{n}^r \subset I\cdot I_{nr-1} + \mathfrak{m}I_{nr},\]
where we used  $\overline{I^r} \subset (x_1, \ldots,
x_n)\overline{I^{r-1}}$. Finally, in $F$, this equation
shows that the indices of nilpotency of the components $F_n$ are
bounded by $e(I)$, as desired. Now we apply
Proposition~\ref{indicesnorm} (and delocalize back to the original
homogeneous ideals). \QED



\subsection*{Cohen-Macaulay normalization}

Expectably, normalization indices are easier to obtain when the
normalization of the ideal is Cohen-Macaulay. The following is
directly derived from the known  characterizations of
Cohen-Macaulayness of Rees algebras of ideals  in terms of
associated graded rings and reduction numbers (\cite{AHT},
\cite{JK95}, \cite{SUV2}).

\begin{Theorem} \label{redintclos} Let $(R,\mathfrak{m})$ be a
Cohen--Macaulay local ring and let $\{ I_n, \ n\geq 0\}$ be a
decreasing multiplicative filtration of ideals, with $I_0=R$,
$I_1=I$, and the property that the corresponding Rees algebra
$B=\sum_{n\geq 0}I_nt^n$ is finite over $A$. Suppose that $\height
I\geq 1$
 and let $J$ be a minimal reduction of $I$. If $B$ is Cohen-Macaulay,
 then
\[{I_{n+1}}= J{I_n}=
I_1{I_n} \  \  \ \  \mbox{{\rm for every \,}} n \geq \ell(I_1)-1,\]
and in particular, $B$ is generated over $R[It]$ by forms  of
 degrees at most $\ell-1=\ell(I_1)-1$,
\[ \sum_{n\geq 0} {I_n}t^n = R[{I_1}t, \ldots,
{I_{\ell-1}}t^{\ell-1}].\]
\end{Theorem}

\smallskip

The proof of Theorem~\ref{redintclos} relies on substituting in any
of the proofs mentioned above (\cite[Theorem 5.1]{AHT}, \cite[Theorem
2.3]{JK95},
\cite[Theorem 3.5]{SUV2}) the $I$-adic filtration $\{I^n\}$ by the filtration
$\{I_n\}$.

\begin{Corollary} \label{normalpowers}
Let $(R, \mathfrak{m})$ be a local analytically unramified normal
Cohen--Macaulay ring and let $I$ be an ideal. If $\overline{A}$ is
Cohen-Macaulay then both indices of normalization $s(I)$ and
$s_0(I)$ are at most $\ell(I)-1$. Moreover, if ${I^n}$ is integrally
closed for $n< \ell(I)$, then $I$ is normal.
\end{Corollary}

\medskip

A case this applies to is that of monomial ideals in a polynomial
ring, since
 the ring $\overline{A}$ is Cohen--Macaulay by Hochster's theorem
(\cite[Theorem 6.3.5]{BH}) (see also \cite{RRV}).

\begin{Example}\rm Let $I=I(\mathcal{C})=(x_1x_2x_5,x_1x_3x_4,x_2x_3x_6,x_4x_5x_6)$
be the edge ideal associated to the
clutter

\bigskip

\bigskip


\bigskip

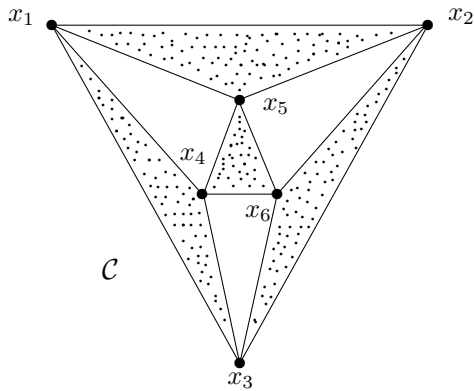
\special{em:linewidth 0.4pt} \unitlength 0.5mm \linethickness{0.4pt}
\begin{picture}(30,100)(-10,30)
\emline{30.00}{130.00}{1}{130.00}{130.00}{2}
\put(80,40){\circle*{2.5}} \put(30,130){\circle*{2.5}}
\put(70,85){\circle*{2.5}} \put(90,85){\circle*{2.5}}
\put(80,110){\circle*{2.5}} \put(130,130){\circle*{2.5}}
\emline{130.00}{130.00}{3}{80.00}{40.00}{4}
\emline{80.00}{40.00}{5}{30.00}{130.00}{6}
\emline{70.00}{85.00}{7}{90.00}{85.00}{8}
\emline{90.00}{85.00}{9}{79.67}{110.00}{10}
\emline{79.67}{110.00}{11}{70.33}{85.00}{12}
\emline{70.33}{85.00}{13}{80.00}{40.00}{14}
\emline{80.00}{40.00}{15}{90.00}{85.00}{16}
\emline{90.00}{85.00}{17}{130.00}{130.50}{18}
\emline{130.00}{130.00}{19}{79.50}{110.00}{20}
\emline{80.00}{110.00}{21}{30.00}{130.00}{22}
\emline{30.00}{130.00}{23}{70.00}{85.00}{24}
\put(67.33,124.00){\circle*{.1}} \put(73.67,124.33){\circle*{.1}}
\put(75.83,124.33){\circle*{.1}} \put(84.67,123.33){\circle*{.1}}
\put(91.33,125.00){\circle*{.1}} \put(90.00,127.67){\circle*{.1}}
\put(89.67,127.67){\circle*{.1}} \put(80.33,119.33){\circle*{.1}}
\put(80.33,119.33){\circle*{.1}} \put(92.33,121.33){\circle*{.1}}
\put(106.33,125.33){\circle*{.1}}
\put(106.33,125.33){\circle*{.1}} \put(99.67,124.67){\circle*{.1}}
\put(99.67,122.33){\circle*{.1}} \put(99.67,122.33){\circle*{.1}}
\put(102.00,127.67){\circle*{.1}}
\put(111.33,127.00){\circle*{.1}}
\put(101.33,127.33){\circle*{.1}} \put(85.67,116.00){\circle*{.1}}
\put(57.33,127.00){\circle*{.1}} \put(46.00,127.33){\circle*{.1}}
\put(50.33,125.00){\circle*{.1}} \put(56.33,122.00){\circle*{.1}}
\put(59.67,128.00){\circle*{.1}} \put(62.00,121.33){\circle*{.1}}
\put(69.00,118.00){\circle*{.1}} \put(78.00,114.67){\circle*{.1}}
\put(71.33,127.33){\circle*{.1}} \put(82.00,127.00){\circle*{.1}}
\put(84.67,119.33){\circle*{.1}} \put(74.00,116.33){\circle*{.1}}
\put(86.33,119.67){\circle*{.1}} \put(62.00,124.67){\circle*{.1}}
\put(116.33,127.33){\circle*{.1}}
\put(103.00,123.00){\circle*{.1}} \put(95.67,123.00){\circle*{.1}}
\put(94.67,127.67){\circle*{.1}} \put(95.33,120.00){\circle*{.1}}
\put(86.67,119.33){\circle*{.1}} \put(90.67,116.67){\circle*{.1}}
\put(81.33,114.33){\circle*{.1}} \put(82.67,126.67){\circle*{.1}}
\put(62.00,125.67){\circle*{.1}} \put(65.67,127.67){\circle*{.1}}
\put(51.00,128.00){\circle*{.1}} \put(40.00,128.33){\circle*{.1}}
\put(38.00,119.00){\circle*{.1}} \put(40.33,115.67){\circle*{.1}}
\put(42.67,112.33){\circle*{.1}} \put(42.67,110.67){\circle*{.1}}
\put(42.67,110.67){\circle*{.1}} \put(46.00,110.00){\circle*{.1}}
\put(43.67,108.67){\circle*{.1}} \put(48.67,105.33){\circle*{.1}}
\put(48.83,105.33){\circle*{.1}} \put(49.00,105.33){\circle*{.1}}
\put(46.00,105.00){\circle*{.1}} \put(47.33,101.67){\circle*{.1}}
\put(47.33,101.67){\circle*{.1}} \put(51.00,101.00){\circle*{.1}}
\put(50.33,97.33){\circle*{.1}} \put(53.00,99.00){\circle*{.1}}
\put(55.00,98.00){\circle*{.1}} \put(52.67,96.33){\circle*{.1}}
\put(52.00,94.00){\circle*{.1}} \put(56.00,94.00){\circle*{.1}}
\put(56.00,94.00){\circle*{.1}} \put(60.33,91.00){\circle*{.1}}
\put(57.67,89.33){\circle*{.1}} \put(58.33,87.00){\circle*{.1}}
\put(60.67,87.00){\circle*{.1}} \put(63.33,87.00){\circle*{.1}}
\put(66.67,86.67){\circle*{.1}} \put(67.67,83.67){\circle*{.1}}
\put(67.67,83.67){\circle*{.1}} \put(63.67,83.00){\circle*{.1}}
\put(59.33,82.33){\circle*{.1}} \put(54.67,90.33){\circle*{.1}}
\put(60.33,84.33){\circle*{.1}} \put(60.33,79.33){\circle*{.1}}
\put(62.67,79.67){\circle*{.1}} \put(65.00,79.67){\circle*{.1}}
\put(67.00,79.67){\circle*{.1}} \put(69.33,80.00){\circle*{.1}}
\put(69.67,76.67){\circle*{.1}} \put(67.00,76.67){\circle*{.1}}
\put(65.00,76.67){\circle*{.1}} \put(65.00,76.67){\circle*{.1}}
\put(62.67,76.67){\circle*{.1}} \put(65.00,72.67){\circle*{.1}}
\put(67.33,70.67){\circle*{.1}} \put(69.67,71.00){\circle*{.1}}
\put(69.00,68.67){\circle*{.1}} \put(69.00,66.67){\circle*{.1}}
\put(71.00,66.00){\circle*{.1}} \put(70.00,64.33){\circle*{.1}}
\put(72.33,63.00){\circle*{.1}} \put(71.67,61.00){\circle*{.1}}
\put(73.67,58.00){\circle*{.1}} \put(73.67,56.67){\circle*{.1}}
\put(75.00,54.33){\circle*{.1}} \put(76.00,51.67){\circle*{.1}}
\put(56.33,87.33){\circle*{.1}} \put(58.67,94.00){\circle*{.1}}
\put(84.33,51.00){\circle*{.1}} \put(84.00,51.67){\circle*{.1}}
\put(87.00,55.67){\circle*{.1}} \put(86.00,58.00){\circle*{.1}}
\put(88.33,57.67){\circle*{.1}} \put(86.67,60.33){\circle*{.1}}
\put(66.67,120.00){\circle*{.1}} \put(54.67,125.00){\circle*{.1}}
\put(106.67,127.67){\circle*{.1}}
\put(108.67,124.00){\circle*{.1}} \put(81.33,123.00){\circle*{.1}}
\put(87.33,125.67){\circle*{.1}} \put(89.00,122.00){\circle*{.1}}
\put(75.00,127.33){\circle*{.1}} \put(82.67,117.00){\circle*{.1}}
\put(92.67,83.33){\circle*{.1}} \put(92.67,79.33){\circle*{.1}}
\put(98.00,83.00){\circle*{.1}} \put(88.00,65.00){\circle*{.1}}
\put(90.00,64.67){\circle*{.1}} \put(92.33,65.67){\circle*{.1}}
\put(94.00,68.33){\circle*{.1}} \put(91.00,69.00){\circle*{.1}}
\put(89.00,70.33){\circle*{.1}} \put(92.00,71.00){\circle*{.1}}
\put(95.00,72.00){\circle*{.1}} \put(92.67,74.00){\circle*{.1}}
\put(90.67,74.33){\circle*{.1}} \put(90.67,76.33){\circle*{.1}}
\put(93.33,76.33){\circle*{.1}} \put(96.33,76.33){\circle*{.1}}
\put(98.33,76.33){\circle*{.1}} \put(99.33,78.33){\circle*{.1}}
\put(96.67,79.00){\circle*{.1}} \put(95.00,80.33){\circle*{.1}}
\put(91.67,81.00){\circle*{.1}} \put(96.00,84.00){\circle*{.1}}
\put(102.00,84.67){\circle*{.1}} \put(98.33,86.00){\circle*{.1}}
\put(95.33,86.67){\circle*{.1}} \put(98.67,88.33){\circle*{.1}}
\put(102.67,88.33){\circle*{.1}} \put(104.67,88.67){\circle*{.1}}
\put(106.00,92.33){\circle*{.1}} \put(102.33,93.67){\circle*{.1}}
\put(100.33,93.67){\circle*{.1}} \put(104.67,95.33){\circle*{.1}}
\put(107.33,96.33){\circle*{.1}} \put(111.00,97.33){\circle*{.1}}
\put(107.67,101.67){\circle*{.1}}
\put(110.67,102.00){\circle*{.1}}
\put(112.33,102.33){\circle*{.1}}
\put(113.67,104.33){\circle*{.1}}
\put(113.33,106.00){\circle*{.1}}
\put(111.33,106.33){\circle*{.1}}
\put(114.00,108.33){\circle*{.1}}
\put(116.33,109.67){\circle*{.1}}
\put(117.00,111.67){\circle*{.1}}
\put(119.00,113.00){\circle*{.1}}
\put(119.67,115.33){\circle*{.1}}
\put(121.67,117.33){\circle*{.1}}
\put(121.67,118.67){\circle*{.1}}
\put(123.67,121.00){\circle*{.1}}
\put(120.33,128.00){\circle*{.1}}
\put(112.67,126.00){\circle*{.1}} \put(86.00,127.33){\circle*{.1}}
\put(70.33,121.33){\circle*{.1}} \put(72.67,120.00){\circle*{.1}}
\put(77.67,121.00){\circle*{.1}} \put(68.67,126.67){\circle*{.1}}
\put(59.33,124.00){\circle*{.1}} \put(68.33,74.00){\circle*{.1}}
\put(62.67,90.00){\circle*{.1}} \put(73.33,87.00){\circle*{.1}}
\put(77.00,89.33){\circle*{.1}} \put(80.00,91.00){\circle*{.1}}
\put(80.33,89.00){\circle*{.1}} \put(84.00,89.00){\circle*{.1}}
\put(84.00,91.33){\circle*{.1}} \put(86.33,87.33){\circle*{.1}}
\put(82.33,95.33){\circle*{.1}} \put(80.33,95.33){\circle*{.1}}
\put(77.33,95.33){\circle*{.1}} \put(79.33,98.33){\circle*{.1}}
\put(80.67,98.67){\circle*{.1}} \put(80.67,100.33){\circle*{.1}}
\put(78.67,100.67){\circle*{.1}} \put(80.00,102.33){\circle*{.1}}
\put(80.00,104.33){\circle*{.1}} \put(80.00,105.67){\circle*{.1}}
\put(95.67,126.00){\circle*{.1}} \put(89.33,61.00){\circle*{.1}}
\put(89.00,67.33){\circle*{.1}} \put(90.33,72.00){\circle*{.1}}
\put(95.67,74.33){\circle*{.1}} \put(100.33,81.33){\circle*{.1}}
\put(93.67,86.00){\circle*{.1}} \put(97.33,90.67){\circle*{.1}}
\put(100.67,90.33){\circle*{.1}} \put(104.33,91.00){\circle*{.1}}
\put(110.00,99.67){\circle*{.1}} \put(105.67,99.00){\circle*{.1}}
\put(103.67,97.00){\circle*{.1}} \put(110.67,104.33){\circle*{.1}}
\put(66.33,67.67){\circle*{.1}} \put(69.67,62.00){\circle*{.1}}
\put(71.33,68.67){\circle*{.1}} \put(71.00,74.00){\circle*{.1}}
\put(65.00,85.00){\circle*{.1}} \put(57.67,84.33){\circle*{.1}}
\put(66.67,82.00){\circle*{.1}} \put(49.33,99.67){\circle*{.1}}
\put(50.67,103.33){\circle*{.1}} \put(46.00,107.67){\circle*{.1}}
\put(60.00,121.00){\circle*{.1}} \put(64.67,119.33){\circle*{.1}}
\put(72.33,115.67){\circle*{.1}} \put(76.33,114.33){\circle*{.1}}
\put(79.33,117.00){\circle*{.1}} \put(76.33,119.33){\circle*{.1}}
\put(74.33,121.33){\circle*{.1}} \put(70.67,124.33){\circle*{.1}}
\put(64.67,123.00){\circle*{.1}} \put(47.33,125.00){\circle*{.1}}
\put(48.33,128.00){\circle*{.1}} \put(43.33,127.33){\circle*{.1}}
\put(54.33,127.33){\circle*{.1}} \put(98.00,127.67){\circle*{.1}}
\put(90.00,119.33){\circle*{.1}} \put(89.00,116.33){\circle*{.1}}
\put(93.00,118.33){\circle*{.1}} \put(98.33,120.67){\circle*{.1}}
\put(106.00,122.33){\circle*{.1}}
\put(103.33,125.67){\circle*{.1}}
\put(108.67,125.67){\circle*{.1}} \put(78.33,87.33){\circle*{.1}}
\put(82.67,87.33){\circle*{.1}} \put(82.00,90.33){\circle*{.1}}
\put(84.33,94.67){\circle*{.1}} \put(81.33,93.00){\circle*{.1}}
\put(81.33,93.00){\circle*{.1}} \put(75.67,91.67){\circle*{.1}}
\put(74.33,89.33){\circle*{.1}} \put(77.67,97.67){\circle*{.1}}
\put(82.67,97.33){\circle*{.1}} \put(81.33,101.67){\circle*{.1}}
\put(79.33,127.33){\circle*{.1}} \put(79.00,124.33){\circle*{.1}}
\put(81.67,125.00){\circle*{.1}} \put(79.00,122.67){\circle*{.1}}
\put(83.33,120.33){\circle*{.1}} \put(53.33,123.33){\circle*{.1}}
\put(80.00,112.67){\circle*{.1}} \put(78.33,92.33){\circle*{.1}}
\put(74.67,90.33){\circle*{.1}} \put(77.67,98.00){\circle*{.1}}
\put(75.67,93.00){\circle*{.1}} \put(78.00,94.67){\circle*{.1}}
\put(75.67,94.33){\circle*{.1}}
\put(22.00,132.00){\makebox(0,0)[cc]{\small $x_1$}}
\put(139.00,132.33){\makebox(0,0)[cc]{\small $x_2$}}
\put(80.33,35.33){\makebox(0,0)[cc]{\small $x_3$}}
\put(89.67,108.33){\makebox(0,0)[cc]{\small $x_5$}}
\put(85.00,80.00){\makebox(0,0)[cc]{\small $x_6$}}
\put(67.67,95.00){\makebox(0,0)[cc]{\small $x_4$}}
\put(45.67,65.00){\makebox(0,0)[cc]{$\mathcal C$}}
\end{picture}

\bigskip

\medskip

Consider the incidence matrix $A$ of this clutter, i.e., the columns 
of $A$ are the exponent vectors of the monomials that generate $I$. Since 
the polyhedron $Q(A)=\{x\vert xA\geq 1;\, x\geq 0\}$ is integral, we 
have the equality $\overline{R[It]}=R_s(I)$, the symbolic Rees
algebra of $I$ (see
\cite[Proposition~3.13]{reesclu}). The ideal $I$ is not normal
because the monomial $x_1x_2\cdots x_6$ is in
$\overline{I^2}\setminus I^2$. 

The minimal primes of $I$ are:
$$
\begin{array}{llll}
\mathfrak{p}_1=(x_1,x_6),&\mathfrak{p}_2=(x_2,x_4),&\mathfrak{p}_3=(x_3,x_5),&
\\
\mathfrak{p}_4=(x_1,x_2,x_5),&\mathfrak{p}_5=(x_1,x_3,x_4),&
\mathfrak{p}_6=(x_2,x_3,x_6),&\mathfrak{p}_7=(x_4,x_5,x_6).
\end{array}
$$
For any $n$, 
\[I^{(n)}=\bigcap_{i=1}^7\mathfrak{p}_i^n.\] A computation
with {\it Macaulay\/} $2$ (\cite{Macaulay2}) gives that
$\overline{I^2}=(I^2,x_1x_2\cdots x_6)$ and that
$\overline{I^3}=I\ \overline{I^2}$. By
Theorem~\ref{redintclos},
$I^{(n)}=II^{(n-1)}$ for $n\geq \ell(I)=4$, where $\ell(I)$ is the
analytic spread of $I$. As a consequence, 
\[\bar{R[It]}= R[It, x_1x_2x_3x_4x_5x_6t^2]. \]

\end{Example}



\begin{Question}{\rm Given the usefulness of
Theorem~\ref{redintclos}, it would be worthwhile to look at the
situation short of Cohen--Macaulayness. For the integral closure of
a standard graded algebra $A$ of dimension $d$,
 it was possible in \cite{emb} to derive
degree bounds assuming only $S_{d-1}$ for $\overline{A}$. Another
issue is to compute the relation type of $A$ in
Theorem~\ref{redintclos}.}\end{Question}



\section{Sally modules and normalization of ideals}\label{sallynormid}


Let $(R, \mathfrak{m})$ be an analytically unramified local ring of
dimension $d$ and $I$ an $\mathfrak{m}$-primary ideal. Let
$\mathcal{F}= \{ I_n, \ n\geq 0\}$ be a decreasing, multiplicative
filtration of ideals, with $I_0=R$, $I_1= I$, with the property that
the corresponding Rees algebra
\[\Rees = R(\mathcal{F})  = \sum_{n\geq 0} I_nt^n \]
is a Noetherian ring. We will examine in detail the case when
$\mathcal{F}$ is a subfiltration of the integral closure filtration
of the powers of $I$, $I_n\subset \overline{I^n}$.

There are several algebraic structures attached to $\mathcal{F}$,
among which we single out  the associated graded ring of
$\mathcal{F}$ and its Sally modules. The first is
\[ \gr_{\mathcal{F}}(R)
= \sum_{n\geq 0} I_n/I_{n+1},\] whose properties are closely linked
to $R(\mathcal{F})$. We note that if $J$ is a minimal reduction of
$I_1$, then $\gr_{\mathcal{F}}(R)$ is a finite generated module over
$\gr_J(R)$, so that it is a semi-standard graded algebra.

To define the Sally module,
 we choose a minimal reduction $J$ of $I$ (if need be, we may
assume that the residue field of $R$ is infinite).
Note that $\Rees$ is a finite extension of the Rees algebra $\Rees_0=
R[Jt]$ of the ideal $J$. The corresponding Sally module $S$ defined by
the
 exact sequence of finitely generated modules over
$\Rees_0$,
\begin{eqnarray}
 0 \rar I\Rees_0 \lar \Rees_{+}[+1] \lar S = \bigoplus_{n=
 1}^{\infty}
I_{n+1}/J^nI \rar 0. \label{sallymod}
\end{eqnarray}

As an $\Rees_0$-module, $S$ is annihilated by an
$\mathfrak{m}$--primary ideal. If $S\neq 0$, $\dim S \leq d$, with
equality if $R$ is Cohen-Macaulay. The
Artinian module
\[S/JtS = \bigoplus_{k\geq 1} I_{k+1}/JI_k \]
gives some control over the number of generators of $\Rees$ as an
$\Rees_0$-module.
 If $S$ is
Cohen-Macaulay, this number is also its multiplicity.
It would, of course, be more useful to obtain bounds for  the
length of $\Rees/(\mathfrak{m}, \Rees_{+})\Rees$, but this requires
lots more.

\medskip

The cohomological properties of $\Rees$, $\gr_{\mathcal{F}}(R)$ and
$S$ become more entwined when $R$ is Cohen-Macaulay. Indeed, under this
condition, the exact sequence (\ref{sallymod}) and
 the
 exact sequences (originally paired in \cite{Hu0}):

\begin{eqnarray}
0 \rar  \Rees_{+}[+1] \lar \Rees \lar {\rm gr}_{\mathcal{F}}(R)
\rar 0\  & & \label{hunekeeq1} \\
& & \nonumber \\
0 \rar \Rees_{+} \lar \Rees \lar R \rar 0,   & & \label{hunekeeq2}
\end{eqnarray}
with the tautological  isomorphism
\begin{eqnarray*}
&&  \Rees_{+}[+1] \cong \Rees_{+}
\end{eqnarray*}
gives a fluid mechanism to pass cohomological information around.

\medskip

\begin{Proposition} Let $(R, \mathfrak{m})$ be a Cohen-Macaulay local
ring of dimension $d$ and $\mathcal{F}$ a filtration as above. Then
\begin{itemize}
\item[{\rm (a)}]
 $\depth \Rees \leq \depth \gr_{\mathcal{F}}(R)+1  $,
 with equality if
 $\gr_{\mathcal{F}}(R)$ is not Cohen-Macaulay.

\item[{\rm (b)}]  $\depth S\leq \depth \gr_{\mathcal{F}}(R) +1 $,
 with equality if
 $\gr_{\mathcal{F}}(R)$ is not Cohen-Macaulay.
\end{itemize}
\end{Proposition}

\demo For (a), see \cite[Theorem 3.25]{bookthree}. For (b), it
follows simply because $I \Rees_0$ is a maximal Cohen--Macaulay
$\Rees_0$-module. \QED

\subsection*{Hilbert functions}

Another connection between $\mathcal{F}$ and $S$ is realized via their
Hilbert functions. Set
\[ H_{\mathcal{F}}(n) = \lambda(R/I_n), \quad H_{S}(n-1)=
\lambda(I_{n}/IJ^{n-1}).\] The associated Poincar\'{e} series
 \[P_{\mathcal{F}}(t) = \frac{f(t)}{(1-t)^{d+1}},\]
 \[P_{S}(t) = \frac{g(t)}{(1-t)^{d}}\]
are related by
\begin{eqnarray*}
 P_{\mathcal{F}}(t) &=&  {\displaystyle\frac{\lambda(R/J)\cdot
 t}{(1-t)^{d+1}}  +
\frac{\lambda(R/I)(1-t)}{(1-t)^{d+1}}
 -P_S(t)}\\
&=& {\displaystyle \frac{\lambda(R/I)
 + \lambda(I/J)\cdot t}{(1-t)^{d+1}} - P_S(t)}.
\end{eqnarray*}

The proof of this fact follows as in  \cite{red} and \cite{Vaz},
replacing the $I$-adic filtration by the filtration $\mathcal{F}$.

\medskip

\begin{Proposition} The $h$-polynomials $f(t)$ and $g(t)$ are related
by
\begin{eqnarray}\label{relh}
f(t) &=&  \lambda(R/I)+ \lambda(I/J)\cdot t - (1-t)g(t).
\end{eqnarray}
In particular, if $f(t)= \sum_{i\geq 0}a_it^i$ and $g(t) =
\sum_{i\geq 1}b_it^i$, then for $i\geq 2$
\[ a_i= b_{i-1}-b_i.\]
\end{Proposition}

\medskip

\begin{Corollary}\label{decreasing} If  $\gr_{\mathcal{F}}(R)$ is
Cohen-Macaulay, then the $h$-vector of $S$ is positive and  non-increasing,
\[ b_i\geq 0, \quad b_1\geq b_2 \geq \cdots \geq 0.\]
In particular, if $b_{k+1}= 0$ for some $k$, then $\Rees$ is generated by its
elements of degree at most $ k$.
\end{Corollary}

\demo That $b_i\geq 0$ follows because $S$ is Cohen--Macaulay, while
the positivity of the $a_i$'s for the same reason and the difference
relation shows that $b_{i-1}\geq b_i$. For the other assertion,
since $S$ is Cohen-Macaulay, $b_k= \lambda(I_{k+1}/JI_k)$. The proof
of this fact is a modification of \cite[See 1.1]{Vaz}, using  the
filtration $\mathcal{F}$ instead of the $I$-adic filtration.
Therefore if $b_k$ vanishes no fresh generators for $I_{k+1}$ are
needed. \QED

\medskip

\begin{Remark}\label{jun18-05}{\rm The equality (\ref{relh}) has several
useful general properties, of which we remark the following. For
$k\geq 2$, one has
\[ f^{(k)}(1) = k g^{(k-1)}(1),\] that is the coefficients $e_i$ of the
Hilbert polynomials of $\gr_{F}(R)$ and $S$ are identical, more
precisely
\[ e_{i+1}(\mathcal{F}) = e_i(S), \quad i\geq 1.\]

Observe that when $\depth \gr_{\mathcal{F}}(R)\geq d-1$, $S$ is
Cohen--Macaulay so its $h$-vector is positive, and therefore all the
$e_i$ along with it (see \cite[Corollary 2]{Marley}). }
\end{Remark}

\medskip

\begin{Corollary}\label{deg4}
If  $\gr_{\mathcal{F}}(R)$ is
Cohen-Macaulay and  $g(t)$ is a polynomial of degree at most $4$,
then
\[ e_{2}(\mathcal{F})
 \geq e_3(\mathcal{F}) \geq e_4(\mathcal{F}) \geq e_5(\mathcal{F}). \]
\end{Corollary}

\demo By our assumption  \[ g(t) = b_1t + b_2t^2 + b_3t^3 +
b_4t^4.\] As
\[ e_{k+1}(\mathcal{F}) = e_{k}(S) = \frac{g^{(k)}(1)}{k!},
\] we have the following equations:
\begin{eqnarray*}
  e_{2}(\mathcal{F})&=& b_1 + 2 b_2 + 3 b_3 + 4 b_4 \\
  e_{3}(\mathcal{F}) &=& b_2 + 3 b_3 + 6 b_4 \\
  e_{4}(\mathcal{F}) &=& b_3 +4 b_4 \\
  e_{5}(\mathcal{F}) &=& b_4
\end{eqnarray*}
Now the assertion follows because $ b_1\geq b_2 \geq b_3 \geq  b_4
\geq 0$, according to Corollary~\ref{decreasing}. \QED

\bigskip

This considerably lowers the possible number of distinct Hilbert
functions for such algebras.

\begin{Remark}
{\rm The assumptions of Corollary~\ref{deg4} are satisfied for
instance if $\dim R\leq 6$ and $\Rees$ is Cohen-Macaulay.
\medskip

These relations provide a fruitful ground for  several questions. Let 
$(R, \mathfrak{m})$ be a local Nagata domain, $I$ an
$\mathfrak{m}$-primary ideal and 
suppose $\mathcal{F}$ is the integral closure filtration of $I$. How
generally does the inequality
\[ \overline{e_1}({I})\geq e_0(I)- \lambda(R/\overline{I})\]
hold? Of course it is true if $\gr_{\mathcal{F}}(R)$ is
Cohen-Macaulay, and possibly if the Sally module $S$ is
Cohen-Macaulay.

}\end{Remark}

\subsection*{Number of generators}

Another application of Sally modules is to obtain a bound for the
number of generators (and the distribution of their degrees)
of $\Rees$ as an $\Rees_0$-algebra. (Sometimes the notation is used to
denote the number of module generators.)
A distinguished feature is the front loading of the new generators in
the Cohen-Macaulay case.

\begin{Theorem} Let $\mathcal{F}$ be a filtration as above.
\begin{itemize}
\item[{\rm (a)}]
 If $\depth \gr_{\mathcal{F}}(R)\geq d-1$,
 the $\Rees_0$-algebra $\Rees$ can be
generated by $e_1(\mathcal{F})$ elements.
\item[{\rm (b)}]
 If $\depth \gr_{\mathcal{F}}(R)= d$,
 the $\Rees_0$-algebra $\Rees$ can be
generated by $e_0(\mathcal{F})$ elements.
\item[{\rm (c)}]
 If $\Rees$ is Cohen-Macaulay, it can be
 generated by $\lambda(\mathcal{F}_1/J) + (d-2)
\lambda(\mathcal{F}_2/\mathcal{F}_1J)$ elements.
\end{itemize}
\end{Theorem}

\demo
 (a) From the relation (\ref{relh}), we
have
\[ e_1(\mathcal{F}) = f'(1) = \lambda(F_1/J) + g(1).\]
From the sequence (\ref{sallymod}) that defines $S$, one has
\[ \nu(\Rees) \leq \nu(F_1/J) + \nu(S)\leq \lambda(F_1/J) + \nu(S).\]
Since $S$ is Cohen-Macaulay, $\nu(S)= e_0(S)=g(1)$, which combine to
give the
 promised assertion.

\medskip

(b) A generating set for $\Rees$ can be obtained from a lift of a
minimal set of generators for $\gr_{\mathcal{F}}(R)$, which is given
by its multiplicity since it is Cohen-Macaulay.

\medskip

(c) Since $\Rees $ is Cohen-Macaulay, its reduction number is $\leq
d-1$. Thus the $h$-polynomial of $\gr_{\mathcal{F}}(R)$ has degree
$\leq d-1$, and consequently the $h$-polynomial of the Sally module
has degree at most $d-2$. As the $h$-vector of
$S$ is decreasing, its multiplicity is at most
$(d-2)\lambda(\mathcal{F}_2/\mathcal{F}_1J)$, and we conclude as in
(a).
 \QED

\begin{Remark}{\rm
A typical application is to the case $d=2$ with $\mathcal{F}$ being
the filtration $F_n= \overline{I^n}$.
}\end{Remark}

\begin{Remark}{\rm There are several relevant issues here. The first, to get bounds for
$e_1(\mathcal{F})$. 
 This
is addressed in \cite{ni1}. For instance, when $R$ is a regular local
ring of characteristic zero, $e_1(\mathcal{F})\leq
\displaystyle{\frac{(d-1)e_0}{2}}$.


}\end{Remark}












\section{One-step normalization of Rees algebras}\label{onestepnorm}

In this section we present one of the rare instances where the
normalization of the Rees ring can be computed in a single step
using an explicit expression as a colon ideal. Our formula applies
to homogeneous ideals that are generated by forms of the same degree
and satisfy some additional assumptions.

The $G_d$ assumption in the theorem means that the minimal number of
generators $\nu(I_{\mathfrak p})$ is at most ${\rm dim}
R_{{\mathfrak p}}$ for every prime ideal ${\mathfrak p}$ containing
$I\/$ with ${\rm dim} R_{{\mathfrak p}}\leq d -1 $. In the proof of
the theorem we use the theory of residual intersections. Let $s$ be an
integer with $s \geq \height \, I\/$. Recall that ${\mathfrak a} \,
\colon I\/$ is a {\it $s$-residual intersection} of $I\/$ if
${\mathfrak a}$ is an $s$-generated $R$-ideal properly contained in
$I\/$ and ${\rm ht} \, {\mathfrak a} \, \colon I \geq s $.

\begin{Theorem}~\label{compint}
Let $k$ be an infinite field, $R=k[\xvec d]$ a positively graded
polynomial ring and $I$ an ideal of height $g$ generated by forms of
degree $\delta$, and set $A=R[It]$. Assume $I$ satisfies $G_d$, ${\rm depth} \ R/I^j
\geq {\rm dim} \ R/I -j+1$ for $1 \leq j \leq d-g$, and $I$ is
normal locally on the punctured spectrum. Let $J$ be a homogeneous
minimal reduction of $I$ and write $\sigma= \sum_{i=1}^{d} {\rm deg}
x_i$. Then
\[
\overline{A}= R[Jt] :_{R[t]} R_{\geq g\delta -\delta -\sigma +1};
\]
in particular $\overline{I}= J :_{R}  R_{\geq g\delta -\delta
-\sigma +1}.$ Furthermore, $I$ is a normal ideal of linear type if
and only if $\delta \leq \frac{\sigma -1}{g-1}$ or $\nu(I) \leq d-1$
\end{Theorem}
\demo We may assume $g \geq 1$ and $d \geq 2$. Notice that $J$ is of
linear type and $R[Jt]$ is Cohen-Macaulay \cite{U94}. Write $\m$ for
the homogeneous maximal ideal of $R$.

If $\ell(I) < d$ then $I=J$. Since $A=R[Jt]$ is Cohen-Macaulay, $I$
is normal on the punctured spectrum, and ${\rm ht} \, \m A \geq 2$,
it follows that $\overline{A}= A= R[Jt]$. On the other hand, $R[Jt]
:_{R[t]} R_{\geq g\delta -\delta -\sigma +1}= R[Jt]$ because ${\rm
ht} (R_{\geq g\delta -\delta -\sigma +1}) R[Jt] \geq {\rm ht} \, \m
R[Jt] \geq 2$ and $R[Jt]$ is Cohen-Macaulay.

Thus we may assume $\ell(I) =d$. Write $\b =R_{\geq g\delta -\delta
-\sigma +1}.$ Notice that $R[Jt] :_{R[t]} \b = R[Jt] :_{{\rm Quot}
(R[t])} \b$ since $d \geq 2$. As the two $R[Jt]$-modules
$\overline{A}$ and $R[Jt] : \b $ satisfy $S_2$ and as locally on the
punctured spectrum of $R$, $I$ is of linear type and normal, it
suffices to prove the equality $\overline{A}_{\m R[Jt]}= (R[Jt] :
\b)_{\m R[Jt]}$.

Let $f_1, \ldots, f_d$ be a general generating set of $J$ consisting
of forms of degree $\delta$ and let $\varphi$ be a minimal
presentation matrix of $f_1, \ldots, f_d$. Notice that the entries
along any column of $\varphi$ are forms of the same degree. One has
 $R[Jt]\cong R[T_1, \ldots, T_d] / I_1( \underline{T} \varphi)$.
Let $K=k(T_1, \ldots, T_d)$ and $B=K[x_1, \ldots, x_d]/I_1(
\underline{T} \varphi)$. Notice that $B_{\m B} \cong R[Jt]_{\m
R[Jt]}$ Since $B$ is a positively graded $K$-algebra with irrelevant
maximal ideal, we conclude that $B$ is a domain of dimension one. In
the ring $K[x_1, \ldots, x_d]$, $I_1( \underline{T} \varphi)= \a :
J$ is a $(d-1)$-residual intersection of $J$, hence of $I$, where
$\a$ is generated by $d-1$ forms of degree $\delta$. From this we
conclude that $\omega_B \cong I^{d-g}/ \a I^{d-g-1} ((d-1)\delta
-\sigma)$ (see \cite{U94}). Thus $a(B)=(g-1)\delta - \sigma$.

Since $B$ is a positively graded $K$-domain, it follows that
$\overline{B}$ is a positively graded $L$-domain for some finite
field extension $L$ of $K$. As ${\rm dim} \, B =1$, $\overline{B}$
is a principal ideal domain, hence $\overline{B} =L[t]$ for some
homogeneous element $t$ of degree $\alpha > 0$. Since the Hilbert
function of $\overline{B}$ as a $B$-module is constant in degrees
divisible by $\alpha$ and zero otherwise, the conductor of $B$ is of
the form $B : \overline{B} = B_{\geq \varepsilon}$ for some
$\varepsilon$, where $\varepsilon={\rm max}\,  \{i \,  | \,
[\overline{B}/B]_i \not=0 \} +1 $. The sequence
\[
0 \longrightarrow B \longrightarrow \overline{B} \longrightarrow
\overline{B}/B \longrightarrow 0
\]
yields an exact sequence
\[
0 \longrightarrow \overline{B}/B \longrightarrow H^{1}_{\m B}(B)
\longrightarrow
 H^{1}_{\m
b}(\overline{B}) \longrightarrow 0.
\]
If $\overline{B}/B \not= 0$, then $a(B) \geq 0$ since
$\overline{B}/B$ is concentrated in non-negative degrees. On the
other hand $a(\overline{B})= -\alpha < 0$. Thus $ \varepsilon =a(B)
+1 = g \delta -\delta -\sigma +1$. Hence $B : \overline{B} = B_{\geq
g\delta -\delta -\sigma +1}$. If on the other hand $\overline{B}/B =
0$, then $a(B)= a(\overline{B})$, hence $g\delta -\delta -\sigma=
-\alpha < 0$. Thus $B_{\geq g\delta -\delta -\sigma +1}=B= B:
\overline{B}$. Therefore in either case
 $B: \overline{B}=B_{\geq g\delta -\delta
-\sigma +1}= \b B $, or equivalently, $\overline{B}= B : \b B$.
Localizing at $\m B$ we conclude that $ \overline{R[It]}_{\m R[Jt]}
= \overline{R[Jt]}_{\m R[Jt]}= (R[Jt] : \b)_{\m R[Jt]}$.
 \QED

\bigskip

\begin{Remark}
{\rm Notice in the previous theorem if $R$ is standard graded then
$R_{\geq g\delta -\delta -\sigma +1}= \m^{g \delta - \delta -\sigma
+ 1 }$}
\end{Remark}


\begin{Remark}
{\rm In the presence of the $G_d$ assumption the depth conditions in
Theorem~\ref{compint} are satisfied, for example if $I$ is perfect
of height 2, $I$ is Gorenstein of height 3, or more generally $I$ is
licci.}
\end{Remark}

\medskip

Another class of ideals satisfying the assumptions of the theorem
are $1$-dimensional ideals:

\begin{Example}{\rm Let $k$ be an infinite field, $R=k[x_1, \ldots,
x_d]$ a standard graded polynomial ring with homogeneous maximal
ideal $\m$, $I$ a $1$-dimensional reduced ideal generated by forms
of degree $\delta$, and $J$ an ideal generated by $d$ general forms
of degree $\delta$ in $I$. Then for every $n$, \[ \overline{I^n}=
J^n : \m^{(d-2)(\delta -1) - 1}.
\]}\end{Example}

\end{document}